\documentclass[12pt]{article}
\usepackage{latexsym,amsmath,amsthm}
\usepackage{amssymb}
\usepackage{graphicx}
\usepackage{appendix}
\usepackage[dvipsnames]{xcolor}
\usepackage{tikz}

\usepackage[colorlinks=true, linkcolor=blue]{hyperref}

\DeclareRobustCommand{\stirlingone}{\genfrac\{\}{0pt}{}}
\DeclareRobustCommand{\stirlingtwo}{\genfrac\{\}{0pt}{}}

\newtheorem{thm}{Theorem}[section]

\newtheorem{lemma}[thm]{Lemma}

\newtheorem{claim}[thm]{Claim}
\newtheorem{obsv}[thm]{Observation}

\theoremstyle{definition}
\newtheorem{algorithm}[thm]{Algorithm}

\theoremstyle{definition}

\theoremstyle{definition}
\newtheorem{question}[thm]{Question}

\theoremstyle{definition}
\newtheorem{defn}[thm]{Definition}

\theoremstyle{definition}

\theoremstyle{definition}

\theoremstyle{definition}

\theoremstyle{definition}
\newtheorem{example}[thm]{Example}

\theoremstyle{remark}

\theoremstyle{remark}

\usepackage[margin=1in]{geometry}

\usepackage{url}
\usepackage{microtype}

\begin{document}

\title{Totally non-negativity of a family of change-of-basis matrices}

\author{David Galvin and Yufei Zhang\thanks{Department of Mathematics,
University of Notre Dame, Notre Dame IN; {\tt dgalvin1}, {\tt yzhang43@nd.edu}. Galvin in part supported by a Simons Collaboration Grant for Mathematicians.}}

\maketitle

\begin{abstract}
Let ${\bf a}=(a_1, a_2, \ldots, a_n)$ and ${\bf e}=(e_1, e_2, \ldots, e_n)$ be real sequences. Denote by $M_{{\bf e}\rightarrow {\bf a}}$ the $(n+1)\times(n+1)$ matrix whose $(m,k)$ entry ($m, k \in \{0,\ldots, n\}$) is the coefficient of the polynomial $(x-a_1)\cdots(x-a_k)$ in the expansion of $(x-e_1)\cdots(x-e_m)$ as a linear combination of the polynomials $1, x-a_1, \ldots, (x-a_1)\cdots(x-a_m)$. By appropriate choice of ${\bf a}$ and ${\bf e}$ the matrix $M_{{\bf e}\rightarrow {\bf a}}$ can encode many familiar doubly-indexed combinatorial sequences, such as binomial coefficients, Stirling numbers of both kinds, Lah numbers and central factorial numbers.   

In all four of these examples, $M_{{\bf e}\rightarrow {\bf a}}$ enjoys the property of total non-negativity --- the determinants of all its square submatrices are non-negative. This leads to a natural question: when, in general, is $M_{{\bf e}\rightarrow {\bf a}}$ totally non-negative?  

Galvin and Pacurar found a simple condition on ${\bf e}$ that characterizes total non-negativity of $M_{{\bf e}\rightarrow {\bf a}}$ when ${\bf a}$ is non-decreasing. Here we fully extend this result. For arbitrary real sequences ${\bf a}$ and ${\bf e}$, we give a condition that can be checked in $O(n^2)$ time that determines whether $M_{{\bf e}\rightarrow {\bf a}}$ is totally non-negative. When $M_{{\bf e}\rightarrow {\bf a}}$ is totally non-negative, we witness this with a planar network whose weights are non-negative and whose path matrix is $M_{{\bf e}\rightarrow {\bf a}}$. When it is not, we witness this with an explicit negative minor.   
\end{abstract}

\medskip


\medskip

\section{Introduction}

Let ${\bf a}=(a_1, a_2, \ldots, a_n)$ and ${\bf e}=(e_1, e_2, \ldots, e_n)$ be real sequences. Each gives rise to a basis for the space of polynomials in one variable ($x$, say) of degree at most $n$, in the following way: set
$$
{\mathcal B}_{\bf a} = \{1, x-a_1, (x-a_1)(x-a_2), \ldots, \prod_{i=1}^n (x-a_i)\}
$$
and define ${\mathcal B}_{\bf e}$ analogously.

Associated with ${\mathcal B}_{\bf a}$ and ${\mathcal B}_{\bf e}$ there is an $(n+1)\times(n+1)$ change-of-basis matrix $M_{{\bf e}\rightarrow {\bf a}}$, whose $(m,k)$ entry ($0 \leq m\leq n$, $0 \leq k \leq n$) is the coefficient of the polynomial $(x-a_1)\cdots(x-a_k)$ in the expansion of $(x-e_1)\cdots(x-e_m)$ as a linear combination of the polynomials in ${\mathcal B}_{\bf a}$. (We adopt the standard convention here that the empty product evaluates to $1$).

Evidently $M_{{\bf e}\rightarrow {\bf a}}$ is lower-triangular, with $1$'s down the main diagonal. There are numerous combinatorial matrices that are instances of $M_{{\bf e}\rightarrow {\bf a}}$, for suitable choices of ${\bf a}$ and ${\bf e}$. We present some examples now (see \cite[Section 1]{GalvinPacurar} for the mostly standard justifications).
\begin{enumerate}
\item Taking $e_i=-1$ and $a_i=0$ for all $i$ yields $M_{{\bf e}\rightarrow {\bf a}}=\left[\binom{m}{k}\right]$, the matrix of binomial coefficients.
\item Taking $e_i=0$ and $a_i=i-1$ for all $i$ yields $M_{{\bf e}\rightarrow {\bf a}}=\left[\stirlingtwo{m}{k}\right]$, the matrix of Stirling numbers of the second kind.
\item Taking $e_i=-(i-1)$ and $a_i=0$ for all $i$ yields $M_{{\bf e}\rightarrow {\bf a}}=\left[\stirlingone{m}{k}\right]$, the matrix of (unsigned) Stirling numbers of the first kind.
\item Taking $e_i=-(i-1)$ and $a_i=i-1$ for all $i$ yields $M_{{\bf e}\rightarrow {\bf a}}=\left[L(m,k)\right]$, the matrix of Lah numbers. 
\item Taking $e_i = i-1-b_i$ and $a_i=i-1$, the $(m,k)$ entry of $M_{{\bf e}\rightarrow {\bf a}}$ is the Rook number $R_{m-k}(B_m)$, where $b_1, b_2, \ldots$ is any non-decreasing sequence of non-negative integers, $B_m$ is the Ferrers board with $m$ columns that has $b_i$ cells in column $i$, and $R_k(B_m)$ is the number of ways of placing $k$ non-attacking rooks on $B_m$.
\end{enumerate}

A matrix is {\em totally non-negative} (TNN) if all its minors (determinants of square sub-matrices) are non-negative. Totally non-negative matrices occur frequently in combinatorics and have been the subject of much investigation. See for example \cite{Brenti2, Brenti3, FominZelevinsky, GascaMicchelli, Skandera} for an overview. In particular, all of the matrices listed above are known to be totally non-negative. Motivated by the first four of the examples above, in \cite{GalvinPacurar} Galvin and Pacurar asked the following question.
\begin{question} \label{quest-ae}
For which pairs of real sequences $({\bf a},{\bf e})$ is $M_{{\bf e}\rightarrow {\bf a}}$ TNN?
\end{question}
A partial answer was given in \cite{GalvinPacurar}. Specifically, a characterization was given of those ${\bf e}$-sequences for which $M_{{\bf e}\rightarrow {\bf a}}$ is TNN, when ${\bf a}$ is weakly increasing. We now briefly describe that characterization, which depends on a notion of {\it restricted growth} (introduced, as far as we are aware, by Arndt \cite{Arndt}).
\begin{defn} \label{def-res-growth}
For weakly increasing ${\bf a}$, say that ${\bf e}$ is a {\em restricted growth sequence} relative to ${\bf a}$ if for each $i \geq 1$ it holds that $e_i \leq a_{f(i)}$, where $f(1)=1$ and for $i \geq 1$
$$
f(i+1) = \left\{
\begin{array}{cl}
f(i) & \mbox{if $e_i < a_{f(i)}$} \\
f(i)+1 & \mbox{if $e_i = a_{f(i)}$}. 
\end{array}
\right.
$$  
\end{defn}
In other words, each $e_i$ is at most a certain cap. The cap for $e_1$ is $a_1$. If $e_1 < a_1$ then the cap for $e_2$ is also $a_1$, while if $e_1=a_1$ then the cap for $e_2$ is $a_2$. In general, the cap for $e_i$ is some $a_{i'}$, and if $e_i<a_{i'}$ then the cap for $e_{i+1}$ is also $a_{i'}$, while if $e_i=a_{i'}$ then the cap for $e_{i+1}$ is $a_{i'+1}$. If ${\bf a}=(0,1,\ldots, n-1, \ldots)$ then a non-negative integer sequence ${\bf e}$ is a restricted growth sequence relative to ${\bf a}$ exactly if it is a restricted growth sequence in the usual sense, that is, one satisfying $e_1=0$ and $e_{i+1} \leq 1+\max_{j=1, \ldots, i} e_j$ for $i \geq 1$.

\begin{thm} \label{thm-mainGP} (\cite[Theorem 1.2]{GalvinPacurar})
Let ${\bf a}$ be a weakly increasing sequence, and let ${\bf e}$ be an arbitrary sequence. Then
\begin{enumerate}
\item \label{mainthmGP-item1} the matrix $M_{{\bf e}\rightarrow {\bf a}}$ is TNN if and only if ${\bf e}$ is a restricted growth sequence relative to ${\bf a}$, and
\item \label{mainthmGP-item2} if ${\bf e}$ is {\em not} a restricted growth sequence relative to ${\bf a}$ then the failure of $M_{{\bf e}\rightarrow {\bf a}}$ to be TNN is witnessed by a negative entry in $M_{{\bf e}\rightarrow {\bf a}}$.
\end{enumerate}
\end{thm}

Galvin and Pacurar left open the major part of Question \ref{quest-ae}, that of determining when $M_{{\bf e}\rightarrow {\bf a}}$ is TNN when no restriction is placed on ${\bf a}$ (or ${\bf e}$). In this note, we completely resolve this question. Our resolution involves extending the notion of restricted growth to the case when ${\bf a}$ is arbitrary, via the following algorithm.
\begin{algorithm} \label{alg-res-growth}
The input is a pair of real sequences ${\bf a}=(a_1, \ldots, a_n)$ and ${\bf e}=(e_1, \ldots, e_n)$. 
\begin{enumerate}
\item Initialize $X=(e_1, \ldots, e_n)$ (an ordered list), and $i=1$.
\item While $i \leq n$, locate the first element, $e'$ say, of $X$ that is greater than or equal to $a_i$.
\begin{enumerate}
\item If there is no such element, then delete the last element of $X$, and increment $i$ by $1$.
\item If $e'=a_i$, then delete $e'$ from $X$, and increment $i$ by $1$.
\item If $e'>a_i$, then STOP and report the current value of $i$.
\end{enumerate}
\item When $i=n+1$, STOP and report the value $n+1$.
\end{enumerate}
\end{algorithm}

\begin{defn} \label{def-res-growth-gen}
For arbitrary real sequences ${\bf a}$ and ${\bf e}$, say that ${\bf e}$ is a {\em restricted growth sequence} relative to ${\bf a}$ if Algorithm \ref{alg-res-growth} (on input ${\bf a}$ and ${\bf e}$) terminates with the report $n+1$. 
\end{defn}

\begin{example} \label{ex-good-quick-alg}
With ${\bf a}=(3,8,7,5,2,7)$ and ${\bf e}=(2,1,3,7,3,4)$ (so $n=6$), Algorithm \ref{alg-res-growth} proceeds as follows:
\begin{itemize}
\item Initially $X=(2,1,3,7,3,4)$, $i=1$ and $a_1=3$. The first element of $X$ which is greater than or equal to $3$ is the third element, which has value $3$, so $e'=3$. 
\item Since $e'=a_1=3$, we update $X$ to $(2,1,7,3,4)$ and set $i=2$.
\item Since there is no $e'$ (in updated $X$) with $e' \geq a_2=8$, we update $X$ to $(2,1,7,3)$ and set $i=3$.
\item Now $e'=7$, so we update $X$ to $(2,1,3)$ and set $i=4$.
\item No element of $X$ is greater than of equal to $a_4=5$, so $X$ becomes $(2,1)$ and $i$ becomes $5$.
\item Now $e'=2$, $X$ updates to $(1)$ and $i$ to $6$.
\item In this final iteration there is no valid $e'$ so $X$ updates to the empty list, $i$ to $7~(=n+1)$, and the algorithm terminates with report $7$. It follows that $(2,1,3,7,3,4)$ is a restricted growth sequence relative to $(3,8,7,5,2,7)$.
\end{itemize}
\end{example}

\begin{example} \label{ex-bad-quick-alg}
With ${\bf a}=(11,8,3,1)$ and ${\bf e}=(10,9,2,1)$ (so $n=4$), Algorithm \ref{alg-res-growth} proceeds as follows:
\begin{itemize}
\item Initially $X=(10,9,2,1)$ and $i=1$; there is no $e'$ with $e' \geq a_1$, so $X$ updates to $(10,9,2)$.
\item At $i=2$ ($a_2=8$) we have $e'=10>8$, so the algorithm terminates with report $2~(<n+1)$. It follows that $(10,9,2,1)$ is not a restricted growth sequence relative to $(11,8,3,1)$. 
\end{itemize}
\end{example}

It is easy to check that if ${\bf a}$ is weakly increasing then the notions of restricted growth given in Definitions \ref{def-res-growth} and \ref{def-res-growth-gen} coincide, so that the following theorem, the main point of this note, is an extension of item 1 of Theorem \ref{thm-mainGP}.
\begin{thm} \label{thm-mainGZ}
Let ${\bf a}$ and ${\bf e}$ be arbitrary real sequences of length $n$. The matrix $M_{{\bf e}\rightarrow {\bf a}}$ is TNN if and only if ${\bf e}$ is a restricted growth sequence relative to ${\bf a}$. 
\end{thm}

\begin{example} \label{ex-mainGZ-good}
(Continuing Example \ref{ex-good-quick-alg}) With ${\bf a}=(3,8,7,5,2,7)$ and ${\bf e}=(2,1,3,7,3,4)$ we have, by Example \ref{ex-good-quick-alg} and Theorem \ref{thm-mainGZ}, that 
$$
M_{{\bf e}\rightarrow {\bf a}} = 
\left[
\begin{array}{ccccccc}
1 & 0 & 0 & 0 & 0 & 0 & 0 \\
1 & 1 & 0 & 0 & 0 & 0 & 0 \\
2 & 8 & 1 & 0 & 0 & 0 & 0 \\
0 & 42 & 12 & 1 & 0 & 0 & 0 \\
0 & 42 & 42 & 10 & 1 & 0 & 0 \\
0 & 210 & 210 & 62 & 9 & 1 & 0 \\
0 & 840 & 840 & 272 & 44 & 12 & 1
\end{array}
\right]
$$
is TNN (and this may be readily verified directly). 
\end{example}

\begin{example} \label{ex-mainGZ-bad}
(Continuing Example \ref{ex-bad-quick-alg}) With ${\bf 
a}=(11,8,3,1)$ and ${\bf e}=(10,9,2,1)$ we have that 
$$
M_{{\bf e}\rightarrow {\bf a}} = 
\left[
\begin{array}{ccccc}
1 & 0 & 0 & 0 & 0 \\
1 & 1 & 0 & 0 & 0  \\
2 & 0 & 1 & 0 & 0  \\
18 & 2 & 1 & 1 & 0  \\
180 & 32 & 4 & 1 & 1
\end{array}
\right]
$$
is not TNN. This may be readily verified directly, for example by noting that the minor corresponding to rows 1 and 2, columns 0 and 1 is $-2$ (recall that our indexing of rows and columns starts at $0$). Note that this is an example of a non-TNN matrix of the form $M_{{\bf e}\rightarrow {\bf a}}$ where the failure to be TNN is not witnessed by a negative entry in the matrix; in the setting of \cite{GalvinPacurar}, where ${\bf a}$ is  weakly increasing, such an example cannot occur.   
\end{example}

Some remarks are in order. Firstly, note that an $(n+1)\times(n+1)$ matrix has $\Theta\left(4^n/\sqrt{n}\right)$ square submatrices and that calculating the determinant of a $k \times k$ matrix requires $\Omega(k^2)$ operations. So the obvious naive algorithm for determining whether such a matrix is TNN requires $\Omega(n^{3/2}4^n)$ operations. Ando \cite[Corollary 2.2]{Ando} established that to test the total non-negativity of a lower-triangular matrix it suffices to test the $\Theta\left(n2^n\right)$ square submatrices whose columns are consecutive and include the first column. Since a proportion $1-o(1)$ of these submatrices have dimension $\Theta(n)$, using Ando's criterion still requires $\Omega(2^n)$ operations. Standing in contrast to this, Theorem \ref{thm-mainGZ} allows one to check whether $M_{{\bf e}\rightarrow {\bf a}}$ is TNN using $O(n^2)$ arithmetic operations. 

Our next remark is that, unsurprisingly, the proof of Theorem \ref{thm-mainGZ} involves the theory of planar networks, and crucially uses the pivot operation of \cite{GalvinPacurar}. We postpone a detailed  discussion of this to Section \ref{sec-results}. For the moment, let us note that our proof passes through an intermediate algorithm that takes as input arbitrary real sequences ${\bf a}$ and ${\bf e}$ of length $n$, and outputs either
\begin{itemize}
\item a planar network with all non-negative weights whose path matrix is $M_{{\bf e}\rightarrow {\bf a}}$, witnessing that $M_{{\bf e}\rightarrow {\bf a}}$ is TNN
\end{itemize}
or
\begin{itemize}
\item a square submatrix of $M_{{\bf e}\rightarrow {\bf a}}$ (with consecutive rows and columns) whose determinant is negative, witnessing that $M_{{\bf e}\rightarrow {\bf a}}$ is not TNN.
\end{itemize}
(See Algorithm \ref{alg-main} and Claim \ref{claim-alg-works}). We will see that when $M_{{\bf e}\rightarrow {\bf a}}$ is not TNN,  the witnessing negative minor can be produced in $O(n^2)$ steps.

Our final remark is that there are various alternate ways of describing the entries of the matrix $M_{{\bf e}\rightarrow {\bf a}}$; these were all observed in \cite{GalvinPacurar}. Denote by $M^{{\bf a}, {\bf e}}(m,k)$ the $(m,k)$ entry of the matrix $M_{{\bf e}\rightarrow {\bf a}}$.
\begin{enumerate}
\item  The numbers $M^{{\bf a}, {\bf e}}(m,k)$ satisfy the recurrence
\begin{equation} \label{Mae-rec}
M^{{\bf a}, {\bf e}}(m,k) = M^{{\bf a}, {\bf e}}(m-1,k-1) + (a_{k+1}-e_m)M^{{\bf a}, {\bf e}}(m-1,k) ~~~\mbox{for $m,k > 0$}
\end{equation}
with initial conditions $M^{{\bf a}, {\bf e}}(0,0) = 1$, $M^{{\bf a}, {\bf e}}(0,k)  = 0$ for $k > 0$ and $M^{{\bf a}, {\bf e}}(m,0) = \prod_{i=1}^m (a_1-e_i)$ for $m > 0$. Various forms of this recurrence have appeared in the literature. As observed in \cite{Gonzales},
with suitable choices of ${\bf a}$ and ${\bf e}$ the recurrence \eqref{Mae-rec} can encode
\begin{itemize}
\item some generalizations of the classical rook numbers \cite{CelesteCorcinoGonzales},
\item the normal order coefficients of the word $(VU)^n$ in the Weyl algebra generated by symbols $V, U$ satisfying $UV - V U = hV^s$ \cite{CelesteCorcinoGonzales},
\item Hsu and Shiue's generalized Stirling numbers \cite{HsuShiue}, 
\item the Jacobi-Stirling numbers (coefficients of the Jacobi differential operator) \cite{CelesteCorcinoGonzales, EverittKwonLittlejohnWellmanYoon},
\end{itemize}
as well as encoding Binomial coefficients, Stirling numbers of both kinds, Lah numbers and rook numbers. Also, with $e_m=0$ for all $m$ and $a_k=(k-1)^2$ for all $k$, \eqref{Mae-rec} encodes the central factorial numbers of Riordan and Carlitz (see e.g. \cite{AlayontKrzywonos}).
\item We have the explicit expression 
\begin{equation*}
M^{{\bf a}, {\bf e}}(m,k) = \sum \prod_{i=1}^{m-k} (a_{s_i-i+1} -e_{s_i})
\end{equation*}
where the sum is over all $S=\{s_1,\ldots,s_{m-k}\} \subseteq \{1,\ldots,m\}$ with $s_1<\ldots<s_{m-k}$.
\item We also have the explicit expression
$$
M^{{\bf a}, {\bf e}}(m,k) = \sum_{\ell=0}^{m-k} (-1)^\ell h_{m-k-\ell}(a_1,\ldots,a_{k+1})s_\ell(e_1,\ldots,e_m)
$$
where $h_{\ell}(x_1,\ldots, x_t)$ is the degree $\ell$ complete symmetric polynomial in $x_1,\ldots, x_t$ and $s_{\ell}(x_1,\ldots, x_t)$ is the degree $\ell$ elementary symmetric polynomial in $x_1,\ldots, x_t$.
\end{enumerate}

In Section \ref{subsec-planar} we review the necessary background material from the theory of planar networks, and describe the pivot operation from \cite{GalvinPacurar}. In Section \ref{subsec-inter-alg} we give the intermediate planar network algorithm. The proof that this algorithm is correct is given in Section \ref{subsec-alg-proof}, while the derivation of Theorem \ref{thm-mainGZ} from all that has preceded appears in Section \ref{subsec-main-proof}.

\section{Planar networks and the intermediate algorithm} \label{sec-results}

\subsection{Setup} \label{subsec-planar}

We require some well-known results from the theory of totally non-negative matrices. (Note that much of what follows is taken from \cite[Section 3]{GalvinPacurar}).

A {\em planar network} $P$ is a directed acyclic planar graph with a subset $\{s_i:i=0,1,2,\ldots\}$ of vertices designated as sources and a subset $\{t_i:i=0,1,2,\ldots\}$ of vertices designated as sinks. We assume throughout that the sources and sinks are arranged around the outer face of the graph, and are ordered cyclically (clockwise) as $\ldots, s_2, s_1, s_0, t_0, t_1, t_2, \ldots$. A {\em weighted} planar network $(P,w)$ is a planar network $P$ together with a function $w:\vec{E}(P)\rightarrow {\mathbb R}$, which we think of as an assignment of weights to the edges of $P$. 

Figure \ref{fig1} shows a particular weighted planar network. The vertices set of this network comprises the sources, the sinks, and the intersection points between the vertical and horizontal lines. All horizontal lines are oriented to the right and all vertical lines are oriented upward. All horizontal edges have weight $1$, while the weights of the vertical edges are given by the $x_{ij}$'s. 

\begin{figure}[ht!]
\begin{center}
\begin{tikzpicture}
\draw [thick, ->] (-6,6) -- (6,6);
\draw [thick, ->] (-6,5) -- (6,5);
\draw [thick, ->] (-6,4) -- (6,4);
\draw [thick, ->] (-6,3) -- (6,3);
\draw [thick, ->] (-6,1) -- (6,1);
\draw [thick, ->] (-6,0) -- (6,0);
\draw [thick, ->] (-5,0) -- (-5,6);
\draw [thick, ->] (-3.5,0) -- (-3.5,5);
\draw [thick, ->] (-2,0) -- (-2,4);
\draw [thick, ->] (-.5,0) -- (-.5,3);
\draw [thick, ->] (5,0) -- (5,1); 
\node [left] at (-6,6) {$s_0$};
\node [left] at (-6,5) {$s_1$};
\node [left] at (-6,4) {$s_2$};
\node [left] at (-6,3) {$s_3$};
\node [left] at (-6,2) {$\vdots$};
\node [left] at (-6,1) {$s_{n-1}$};
\node [left] at (-6,0) {$s_n$};
\node [right] at (6,6) {$t_0$};
\node [right] at (6,5) {$t_1$};
\node [right] at (6,4) {$t_2$};
\node [right] at (6,3) {$t_3$};
\node [right] at (6,2) {$\vdots$};
\node [right] at (6,1) {$t_{n-1}$};
\node [right] at (6,0) {$t_n$};
\node [right] at (-5,5.5) {$\!x_{11}$};
\node [right] at (-5,4.5) {$\!x_{21}$};
\node [right] at (-5,3.5) {$\!x_{31}$};
\node [right] at (-5,.5) {$\!x_{n1}$};
\node [right] at (-3.5,4.5) {$\!x_{22}$};
\node [right] at (-3.5,3.5) {$\!x_{32}$};
\node [right] at (-3.5,.5) {$\!x_{n2}$};
\node [right] at (-2,3.5) {$\!x_{33}$};
\node [right] at (-2,.5) {$\!x_{n3}$};
\node [right] at (-.5,.5) {$\!x_{n4}$};
\node [right] at (5,.5) {$\!x_{nn}$};
\node at (1,2) {$\ddots$};
\end{tikzpicture}
\caption{A weighted planar network. The labels on the vertical edges indicate the value of $w$ for those edges; $w$ takes the value $1$ for all horizontal edges.} \label{fig1}
\end{center}
\end{figure}

The weighted planar network shown in Figure \ref{fig1} is the only one that will be used in our proofs (with various choices of $w$). In figures we will use the compact array representation shown in Figure \ref{array-of-weights}. 

\begin{figure}[ht!]
\begin{equation*}
\begin{array}{ccccccc}
x_{11} & & & & & & \\
x_{21} & x_{22} & & & & & \\
x_{31} & x_{32} & x_{33} & & & & \\
x_{41} & x_{42} & x_{43} & x_{44} & & & \\
x_{51} & x_{52} & x_{53} & x_{54} & x_{55} &   & \\
\vdots & \vdots & \vdots & \vdots &        & \ddots & \\
x_{n1} & x_{n2} & x_{n3} & x_{n4} & \cdots & x_{n(n-1)} & x_{nn}.
\end{array}
\end{equation*}
\caption{A compact array representation of the weighted planar network of Figure \ref{fig1}. We will use the shorthand $A$ for this array.} 
\label{array-of-weights}
\end{figure}

We will say that the edge of the planar network that is given weight $x_{mk}$ in Figure \ref{array-of-weights} is the {\it $[m,k]$ edge} of the array (using square brackets to distinguish this from an entry in a matrix). For each $k=1, \ldots, n$, the collection of $[m,k]$ edges as $m$ ranges from $k$ to $n$ is the {\em $k$-th column} of the planar network (so, for example, the set of edges with weights $x_{33}, \ldots, x_{n3}$ is the third column of the array in Figure \ref{array-of-weights}).

The {\em path matrix} $M(P,w)$ of a generic weighted planar network $(P,w)$ is the $(n+1)\times(n+1)$ matrix whose $(m,k)$ entry (with rows and columns indexed by $\{0,1,2,\ldots,n\}$) is the sum of the weights of all the directed paths from $s_m$ to $t_k$, where the weight of any such path is the product, over all edges traversed, of the weight of the edge. For example, the $(3,1)$ entry of the path matrix of the weighted planar network shown in Figure \ref{fig1} is $x_{31}x_{21}+ x_{31}x_{22}+ x_{32}x_{22}$. Notice that the path matrix of the weighted planar network shown in Figure \ref{fig1} is lower-triangular with $1$'s down the main diagonal. The following \cite{Lindstrom} (see also, for example, \cite[Chapter 29]{AignerZiegler}, \cite{Skandera}) is a standard result in the theory of totally non-negative matrices.
\begin{lemma} (Lindstr\"om's Lemma) \label{lem-tp}
If $w \geq 0$ (i.e., all weights are non-negative) then the matrix $M(P,w)$ is TNN. 

Moreover, the minor corresponding to selecting the rows indexed by $I$ and columns indexed by $J$ (with indexing of rows and columns starting from $0$) equals the sum of the weights of all the collections of $|I|$ vertex disjoint paths from the sources $\{s_i: i \in I\}$ to the sinks $\{t_j: j \in J\}$, where the weight of a collection of paths is the product of the weights of the individual paths in the collection.
\end{lemma}

As mentioned earlier, all the weighted planar networks that appear in our proofs have the form shown in Figure \ref{fig1}. Such networks are completely specified by the array of weights shown in Figure \ref{array-of-weights}. Denoting that array by $A$, we use the shorthand $M(A)$ for $M(P,w)$ (where $(P,w)$ is the weighted planar network shown in Figure \ref{fig1}). Abusing notation slightly, we will refer to $M(A)$ as the path matrix of the array $A$.

All of the arrays that appear in our proofs will be generated from real sequences ${\bf a}=(a_1, \ldots, a_n)$ and ${\bf e}=(e_1, \ldots, e_n)$. Figure \ref{orig-network} shows the generic array $A({\bf a},{\bf e})$ generated from ${\bf a}$ and ${\bf e}$, while Figure \ref{good-network-1-ex} shows a specific instance.

\begin{figure}[ht!]
\begin{equation*}
\begin{array}{ccccccc}
a_1-e_1 & & & & & & \\
a_1-e_2 & a_2-e_1 & & & & & \\
a_1-e_3 & a_2-e_2 & a_3-e_1 & & & & \\
a_1-e_4 & a_2-e_3 & a_3-e_2 & a_4-e_1 & & & \\
a_1-e_5 & a_2-e_4 & a_3-e_3 & a_4-e_2 & a_5-e_1 & & \\
\vdots  & \vdots  & \vdots  & \vdots  &  & \ddots & \\
a_1-e_n & a_2-e_{n-1} & a_3-e_{n-2} & a_4-e_{n-3} & \cdots & a_{n-1}-e_2 & a_n-e_1 
\end{array}
\end{equation*}
\caption{The array $A({\bf a},{\bf e})$.} \label{orig-network}
\end{figure}

\begin{figure}[ht!]
\begin{equation*}
\begin{array}{clc}
\begin{array}{ccccccc}
3-2 & & & & & & \\
3-1 & 8-2 & & & & & \\
3-3 & 8-1 & 7-2 & & & & \\
3-7 & 8-3 & 7-1 & 5-2 & & & \\
3-3 & 8-7 & 7-3 & 5-1 & 2-2 & & \\
3-4 & 8-3 & 7-7 & 5-3 & 2-1 & 7-2 
\end{array}
&
\mbox{or}~~~
&
\begin{array}{rrrrrrr}
1 & & & & & & \\
2 & 6 & & & & & \\
0 & 7 & 5 & & & & \\
-4 & 5 & 6 & 3 & & & \\
0 & 1 & 4 & 4 & 0 & & \\
-1 & 5 & 0 & 2 & 1 & 5 
\end{array}
\end{array}
\end{equation*}
\caption{The array $A({\bf a},{\bf e})$ in the specific case ${\bf a}=(3,8,7,5,2,7)$, ${\bf e}=(2,1,3,7,3,4)$ (Example \ref{ex-good-quick-alg}).} \label{good-network-1-ex}
\end{figure}

We need the following lemma (\cite[Lemma 3.5]{GalvinPacurar}):
\begin{lemma} \label{lemma2-initial_network_construction}
For arbitrary ${\bf a}$ and ${\bf e}$, we have $M(A({\bf a},{\bf e})) = M_{{\bf e}\rightarrow {\bf a}}$.
\end{lemma}

If it happens that all weights in $A({\bf a},{\bf e})$ are non-negative, then the total non-negativity of $M_{{\bf e}\rightarrow {\bf a}}$ follows immediately from Lemmas \ref{lem-tp} and Lemma \ref{lemma2-initial_network_construction}. If $A({\bf a},{\bf e})$ has some negative weights, then we will give a procedure that modifies the array by iteratively replacing negative weights with $0$ weights, while not changing the path matrix of the array. We will show that this procedure succeeds in replacing all negative weights with $0$ weights exactly when $M_{{\bf e}\rightarrow {\bf a}}$ is totally non-negative. To describe the procedure we need (\cite[Lemma 3.9]{GalvinPacurar}), whose statement requires a little more notation. For an array $A$ of the kind shown in Figure \ref{array-of-weights}, we refer to the triangle consisting of the $[m+\ell_1,k+\ell_2]$ edges for $0 \leq \ell_1 \leq n-m$ and $0 \leq \ell_2 \leq \ell_1$ as the triangle {\em headed} at the $[m,k]$ edge (see Figure \ref{fig2}).

\begin{figure}[ht!]
\begin{center}
$$
\begin{array}{ccccccccc}
 & [1,1] & & & & & & & \\
& [2,1] & [2,2] & & & & & & \\
 & [3,1] & {\bf [3,2]} & [3,3] & & & & & \\
 & [4,1] & {\bf [4,2]} & {\bf [4,3]} & [4,4] & & & & \\
 & [5,1] & {\bf [5,2]} & {\bf [5,3]} & {\bf [5,4]} & [5,5] & & & \\
 & \vdots & \vdots &  &  & \ddots & \ddots  & & \\
 & [n,1] & {\bf [n,2]} & {\bf [n,3]} & {\bf [n,4]} & {\bf \cdots} & {\bf [n,n-1]} & [n,n] &
\end{array}
$$
\caption{The triangle headed at the $[3,2]$ position ({\bf bolded} entries).} \label{fig2}
\end{center}
\end{figure}

Suppose that the triangle headed at the $[m,k]$ edge of $A$ is an instance of $A({\bf f},{\bf g})$ for some sequences ${\bf f}=(f_1, \ldots, f_{n-m+1})$, ${\bf g}=(g_1, \ldots, g_{n-m+1})$. Denote by $A^{\rm p}(m,k)$ the array obtained from $A$ by the following operations:
\begin{itemize}
\item Leave unchanged all weights of $A$ that are on edges not in the the triangle headed at the $[m,k]$ edge.
\item Also leave unchanged the weight ($f_1-g_1$) on the $[m,k]$ edge. Replace the weights on the $[m+1,k]$ and $[m+1,k+1]$ edges ($f_1-g_2$ and $f_2-g_1$, respectively), with $f_1-g_1$ and $f_2-g_2$, respectively.
\item More generally, replace the weights on the $[m+\ell,k], [m+\ell,k+1], [m+\ell,k+2], \ldots, [m+\ell,k+\ell]$ edges ($f_1-g_{\ell+1}, f_2-g_\ell, \ldots, f_{\ell+1}-g_1$, respectively), with $f_1-g_1, f_2--g_{\ell+1}, \ldots, f_{\ell+1}-g_2$, respectively, for $1 \leq \ell \leq n-m$.
\end{itemize}
We refer to $A^{\rm p}(m,k)$ as the array obtained from $A$ by {\it pivoting} on the $[m,k]$ edge. Note that in each row of the triangle headed at the $[m,k]$ edge, $A^{\rm p}({\bf f}, {\bf g})$ is obtained from $A$ by moving the $-g_1$'s from the last position in the row to the first, and then shifting all other $-g_j$'s in the row one place to the right. Note also that if $f_1=g_1$, then after pivoting the weights on the $[m',k]$ edges, for all $m'\geq m$, are $0$. (See Figure \ref{good-network-2-ex} for an illustration.)

\begin{figure}[ht!]
\begin{equation*}
\begin{array}{clc}
\begin{array}{ccccccc}
3-2 & & & & & & \\
3-1 & 8-2 & & & & & \\
\underline{{\bf 3-3}} & 8-1 & 7-2 & & & & \\
\underline{3-3} & \underline{8-7} & 7-1 & 5-2 & & & \\
\underline{3-3} & \underline{8-3} & \underline{7-7} & 5-1 & 2-2 & & \\
\underline{3-3} & \underline{8-4} & \underline{7-3} & \underline{5-7} & 2-1 & 7-2 
\end{array}
&
\mbox{or}~~~
&
\begin{array}{rrrrrrr}
1 & & & & & & \\
2 & 6 & & & & & \\
0 & 7 & 5 & & & & \\
0 & 1 & 6 & 3 & & & \\
0 & 5 & 0 & 4 & 0 & & \\
0 & 4 & 4 & -2 & 1 & 5 
\end{array}
\end{array}
\end{equation*}
\caption{The array from Figure \ref{good-network-1-ex} (Example \ref{ex-good-quick-alg}), after pivoting on the $[3,1]$ edge ({\bf bolded}). The edges of the triangle headed at the $[3,1]$ edge (the only edges whose weights might change by pivoting) are \underline{underlined}.} \label{good-network-2-ex}
\end{figure}

\begin{figure}[ht!]
\begin{equation*}
\begin{array}{clc}
\begin{array}{ccccccc}
3-2 & & & & & & \\
3-1 & 8-2 & & & & & \\
3-3 & 8-1 & 7-2 & & & & \\
3-3 & 8-7 & 7-1 & 5-2 & & & \\
3-3 & 8-3 & \underline{{\bf 7-7}} & 5-1 & 2-2 & & \\
3-3 & 8-4 & \underline{7-7} & \underline{5-3} & 2-1 & 7-2 
\end{array}
&
\mbox{or}~~~
&
\begin{array}{ccccccc}
1 & & & & & & \\
2 & 6 & & & & & \\
0 & 7 & 5 & & & & \\
0 & 1 & 6 & 3 & & & \\
0 & 5 & 0 & 4 & 0 & & \\
0 & 4 & 0 & 2 & 1 & 5 
\end{array}
\end{array}
\end{equation*}
\caption{The array from Figure \ref{good-network-2-ex} (Example \ref{ex-good-quick-alg}), after pivoting on the $0$ at the $[5,3]$ position ({\bf bolded}).} \label{good-network-3-ex}
\end{figure}

\begin{lemma} \label{lemma-pivot}
(\cite[Lemma 3.9]{GalvinPacurar}) For an array $A$ of the kind shown in Figure \ref{array-of-weights}, suppose that the triangle headed at the $[m,k]$ edge of $A$ is an instance of $A({\bf f},{\bf g})$ for some sequence ${\bf f}=(f_1, \ldots, f_{n-m+1})$, ${\bf g}=(g_1, \ldots, g_{n-m+1})$. 

If $f_1=g_1$ then $M(A^{\rm p}(m,k)) = M(A)$ (that is, the path matrix is unchanged after pivoting on the $[m,k]$ edge, as long as the entry on the $[m,k]$ edge is $0$).
\end{lemma}

It should be noted that when Lemma \ref{lemma-pivot} was stated in \cite{GalvinPacurar}, some assumptions were made about the weights in $A$ on edges other than those edges in the triangle headed at the $[m,k]$ edge. However, those assumptions were not actually used in the proof, so the proof goes through unchanged in the setting of Lemma \ref{lemma-pivot}.

\subsection{The intermediate algorithm} \label{subsec-inter-alg}

We are now in a position to describe the intermediate planar network algorithm alluded to in the introduction. 

\begin{algorithm} \label{alg-main}
The input is a pair of real sequences ${\bf a}=(a_1, \ldots, a_n)$ and ${\bf e}=(e_1, \ldots, e_n)$. The algorithm proceeds by modifying the weighted planar network $A({\bf a},{\bf e})$ (see Figure \ref{orig-network}).
\begin{itemize}
\item Initially, set $A=A({\bf a},{\bf e})$ and $k=1$.
\item As long as $k \leq n$, scan the weights in the $k$th column of $A$, from the $[k,k]$ edge down.
\begin{itemize}
\item If no zero or negative weights are encountered in the column, leave $A$ unchanged and increase $k$ by $1$.
\item If the first non-positive weight that is encountered is negative then STOP, output the current $A$ and report that $M_{{\bf e}\rightarrow {\bf a}}$ is not TNN.
\item If the first non-positive weight that is encountered is zero, say on the $[m,k]$ edge of the array $A$:
\begin{itemize}
\item if the triangle headed at the $[m,k]$ edge of $A$ is an instance of $A({\bf f},{\bf g})$ for some sequences ${\bf f}=(f_1, \ldots, f_{n-m+1})$, ${\bf g}=(g_1, \ldots, g_{n-m+1})$, then update $A$ by pivoting on the $[m,k]$ position --- that is, replace $A$ by $A^{\rm p}(m,k)$ --- and increase $k$ by $1$. (See Figure \ref{orig-network} for an illustration of an instance of $A({\bf f},{\bf g})$.)
\item if not, then STOP and report that the algorithm is inconclusive. (We will see in the analysis that this clause of the algorithm will never be invoked). 
\end{itemize}
\end{itemize}
\item If $k=n+1$, then STOP, output the current $A$ and report that $M_{{\bf e}\rightarrow {\bf a}}$ is TNN.
\end{itemize}
\end{algorithm}

For the input ${\bf a}=(3,8,7,5,2,7)$, ${\bf e}=(2,1,3,7,3,4)$ (Example \ref{ex-good-quick-alg}), the evolution of the planar network in Algorithm \ref{alg-main} is illustrated in Figures \ref{good-network-1-ex}, \ref{good-network-2-ex} and \ref{good-network-3-ex}. In this case the algorithm reports that $M_{{\bf e}\rightarrow {\bf a}}$ is TNN, and outputs the final network that is illustrated in Figure \ref{good-network-3-ex}. For the input ${\bf a}=(11,8,3,1)$ and ${\bf e}=(10,9,2,1)$ (Example \ref{ex-bad-quick-alg}), Algorithm \ref{alg-main} terminates at $k=2$ (without having done a pivot), and in this case the algorithm reports that $M_{{\bf e}\rightarrow {\bf a}}$ is not TNN. The final network that the algorithm outputs is illustrated in Figure \ref{bad-network-ex}.

\begin{figure}[ht!]
\begin{equation*}
\begin{array}{clc}
\begin{array}{cccc}
11-10 & & &  \\
11-9 & 8-10 & &  \\
11-2 & 8-9 & 3-10 &  \\
11-1 & 8-2 & 3-9 & 1-10 
\end{array}
&
\mbox{or}~~~
&
\begin{array}{rrrr}
1 & & & \\
2 & -2 & & \\
9 & -1 & -7 &  \\
10 & 7 & -6 & -9  
\end{array}
\end{array}
\end{equation*}
\caption{The final network in Algorithm \ref{alg-main} when ${\bf a}=(11,8,3,1)$ and ${\bf e}=(10,9,2,1)$ (Example \ref{ex-bad-quick-alg}).} \label{bad-network-ex}
\end{figure}

\begin{claim} \label{claim-alg-works}
For every input ${\bf a}$, ${\bf e}$:
\begin{enumerate}
\item Algorithm \ref{alg-main} terminates after at most $O(n^3)$ operations.
\item It never terminates with an inconclusive report.
\item If it terminates with the report that $M_{{\bf e}\rightarrow {\bf a}}$ is TNN, then the outputted weighted planar network $A$ has all non-negative entries and has path matrix $M_{{\bf e}\rightarrow {\bf a}}$ (and so serves as a witness that $M_{{\bf e}\rightarrow {\bf a}}$ is TNN).
\item If it terminates with the report that $M_{{\bf e}\rightarrow {\bf a}}$ is not TNN, then there is an $O(n^2)$ algorithm that identifies, in the outputted weighted planar network $A$, a collection $S$ of consecutive sources and a collection $T$ of consecutive sinks with $|S|=|T|$, such that all collections of vertex disjoint paths from $S$ to $N$ have non-positive weight, and at least one has negative weight (which serves, via Lemma \ref{lem-tp}, to witness that $M_{{\bf e}\rightarrow {\bf a}}$ is not TNN).
\end{enumerate}     
\end{claim}

\section{Proofs of Claim \ref{claim-alg-works} and Theorem \ref{thm-mainGZ}}

\subsection{Proof of Claim \ref{claim-alg-works}} \label{subsec-alg-proof}

\subsubsection{Item 1}

Constructing the initial weighted planar network $A$ requires $O(n^2)$ operations. The algorithm requires scanning $O(n^2)$ entries. At most $O(n)$ of these scans require performing a pivot, and each pivot requires $O(n^2)$ recalculations of weights. So in total, the algorithm requires $O(n^3)$ operations.

\subsubsection{Item 2}

This is an immediate corollary of the following straightforward observation:
\begin{obsv} \label{obsv-pivot-good}
Suppose that at some point in the algorithm the triangle headed at the $[k,k]$ edge is of the form $A({\bf f},{\bf g})$ for some sequences ${\bf f}=(f_k, \ldots, f_n), {\bf g}=(g_1, \ldots, g_{n-k+1})$. If we then pivot on the $[m,k]$ edge for some $m \geq k$ (so necessarily $f_k=g_{m-k+1}$), then after the pivot
\begin{itemize}
\item all the entries in positions $[m',k]$ for $m' \geq m$ are $0$ (more precisely, they are all $f_k-g_{m-k+1}$), and
\item the triangle headed at the $[k+1,k+1]$ position is of the form $A({\bf f}',{\bf g}')$ where ${\bf f}'=(f_{k+1}, \ldots, f_n)$ and ${\bf g}'=(g_1, \ldots, g_{m-k}, \widehat{g_{m-k+1}}, g_{m-k+2}, \ldots, g_{n-k+1})$. (The $\widehat{g_{m-k+1}}$ here indicates the removal of that entry from the sequence).
\end{itemize}
\end{obsv}
See Figures \ref{good-network-2-ex} and \ref{good-network-3-ex} for illustrations of Observation \ref{obsv-pivot-good}. 

\subsubsection{Item 3}

Here Observation \ref{obsv-pivot-good} is also a key ingredient. The proof of the statement in item 3 is by induction (on the column $k$ being scanned).

By Lemma \ref{lemma2-initial_network_construction}, the initial network has path matrix $M_{{\bf e}\rightarrow {\bf a}}$, and the triangle headed at the $[1,1]$ edge is of the form $A({\bf f},{\bf g})$ (specifically with ${\bf f}={\bf a}$ and ${\bf g}={\bf e}$).

Suppose that we are at the point in the algorithm where we are scanning the $k$th column. Let us assume (by induction) that in the current weighted planar network $A$, the triangle headed at the $[k,k]$ edge is of the form $A({\bf f},{\bf g})$, and that all the entries in the first $k-1$ columns are non-negative.

If column $k$ has only positive entries, then we move on to column $k+1$ without changing $A$; so as we scan column $k+1$, in the current weighted planar network $A$, the triangle headed at the $[k+1,k+1]$ edge is of the form $A({\bf f}',{\bf g}')$ (with ${\bf f}'$ obtained from ${\bf f}$ by removing the first element of ${\bf f}$, and ${\bf g}'$ obtained from ${\bf g}$ by removing the last element of ${\bf g}$). Also, all the entries in the first $k$ columns are non-negative.

If column $k$ has its first zero at position $[m,k]$, for some $m \geq k$, then, by the observation above, after pivoting the $k$th column has all its entries non-negative, and (in the updated network) the triangle headed at the $[k+1,k+1]$ edge is of the form $A({\bf f}',{\bf g}')$. Also (crucially), by Lemma \ref{lemma-pivot}, the new network still has path matrix $M_{{\bf e}\rightarrow {\bf a}}$.

So by induction (on $k$), when we get to the point where $k=n+1$ (which we do reach, by assumption), the final weighted planar network $A$ has all non-negative entries and has path matrix $M_{{\bf e}\rightarrow {\bf a}}$. 

\subsubsection{Item 4}

Suppose that the first negative weight that is encountered in running the algorithm is at the $[m,k_0]$ edge. 

Mark the $[m,k_0]$ edge. Consider the path $p_m$, from source $s_m$ to sink $t_k$, where $k=(k_0)=k_0-1$, constructed as follows: start by going horizontally from $s_m$ to the lower vertex of the marked edge, then go vertically until the path currently being constructed intersects with the unique path from $s_{k_0-1}$ to $t_{k_0-1}$, and then go horizontally to $t_{k_0-1}$. This path only uses edges with non-zero weight (the $[m',k_0]$ edges, for $m' < m$, have positive weight), and in fact uses only edges with positive weight, except for the $[m,k_0]$ edge, which has negative weight.

Next, let $k_1$ be maximal subject to the conditions that firstly $k_1 < k_0$, and secondly that the $[m-1,k_1]$ edge has positive weight (assuming such a $k_1$ exists --- if not, the process terminates, as described in more detail in a moment), and mark the $[m-1,k_1]$ edge. Note that there is a path,  $p_{m-1}$, from source $s_{m-1}$ to sink $t_{k-1}$ 
--- horizontal to the marked edge, then vertical to the unique path from $s_{k-1}$ to $t_{k-1}$ (the vertical stretch of the network goes at least as high as $k-1$, since it is strictly to the left of a vertical stretch that goes at least as far as the unique path from $s_k$ to $t_k$), then horizontal. This path is vertex disjoint from $p_k$. Also, this path only uses edges with positive weight --- that the marked edge has positive weight forces all edges above it to have positive weight.

Repeat --- mark (if it exists) the $[m-2,k_2]$ edge, where $k_2 < k_1$ is maximal subject to the that edge having positive weight. Note that there is a path, $p_{m-2}$  say, from source $s_{m-2}$ to sink $t_{k-2}$, that is vertex-disjoint from $p_m, p_{m-1}$, and again only uses edges with positive weight.

Repeat, until no markable edge is found. Suppose that the last marked edge is the $[m-\ell,k_\ell]$ edge, where $\ell \geq 0$ (so $\ell+1$ edges have been marked).

There is a collection of vertex disjoint paths in the network $A$ from $\{s_{m-\ell}, \ldots, s_m\}$ to $\{t_{k-\ell}, \ldots, t_k\}$ --- namely $\{p_{m-\ell}, \ldots, p_m\}$ --- and it has negative weight, since all edges on it have positive weight, except the first marked edge at position $[m,k_0]$, which is negative. Note also that the construction of this collection can be done by looking at the weight on each edge of the network only once, and so in time $O(n^2)$.

There may be more such collections of vertex disjoint paths from $\{s_{m-\ell}, \ldots, s_m\}$ to $\{t_{k-\ell}, \ldots, t_k\}$, that avoid edges with weight zero, but we assert that all of these will use the  $[m,k_0]$ edge and no other edge with negative weight. If this were true, then all such collections of paths would have negative weight, and so the sum of the weights of all the collections of vertex disjoint paths from $\{s_{m-\ell}, \ldots, s_m\}$ to $\{t_{k-\ell}, \ldots, t_k\}$ would be negative. From this, it would follow (by Lindstr\"om's lemma, Lemma \ref{lem-tp}) that the minor of the path matrix of $A$ corresponding to rows $m-\ell$ through $m$, columns $k-\ell$ through $k$, would be negative. Since that path matrix is $M_{{\bf e}\rightarrow {\bf a}}$ (this follows from the discussion of item 3), item 4 would follow.

\begin{example}
For the final network that comes from Algorithm \ref{alg-main} with the input ${\bf a}=(11,8,3,1)$ and ${\bf e}=(10,9,2,1)$ (Example \ref{ex-bad-quick-alg}, see Figure \ref{bad-network-ex}), we mark the $[2,2]$ and $[1,1]$ edges, and the identified set of sources is $\{2,1\}$ while the sinks are $\{0,1\}$; this recovers the negative minor pointed out in the introduction.   
\end{example}

\begin{example}
For a more substantial example, consider the weighted planar network whose sign-pattern in shown in Figure \ref{fig-sign-pat}.
\begin{figure}[ht!]
\begin{equation*}
\begin{array}{ccccccc}
+ &  &  &  &  &  &  \\
+ & + &  &  &  &  &  \\
+ & + & + &  &  &  &  \\
+ & + & + & + &  &  &  \\
+ & 0 & + & + & + &  &  \\
+ & 0 & + & + & 0 & + &  \\
+ & 0 & + & 0 & 0 & - & \star  \end{array}
\end{equation*}
\caption{An example of the possible sign pattern on a final planar network from Algorithm \ref{alg-main}, when the algorithm terminates with the report that $M_{{\bf e}\rightarrow {\bf a}}$ in not TNN.} \label{fig-sign-pat}
\end{figure}
In this case the marked edges are those at positions $[7,6]$, $[6,4]$, $[5,3]$, $[4,2]$ and $[3,1]$, the identified set of sources is $\{3,4,5,6,7\}$ and the identified set of sinks is $\{1,2,3,4,5\}$. There is more than one collection of vertex disjoint paths from these sources to these sinks that avoids edges with weight $0$, but all of them use exactly one edge with negative weight (the one at position $[7,6]$), so all have negative weight. 
\end{example}

So it remains to prove the assertion that all collections of vertex disjoint path systems from $\{s_m, \ldots, s_{m-\ell}\}$ to $\{t_k, \ldots, t_{k-\ell}\}$, that avoid edges with weight zero, use the $[m,k_0]$ edge and no other edge with negative weight.

For the proof, first observe that any collection of vertex disjoint paths from $\{s_m, \ldots, s_{m-\ell}\}$ to $\{t_k, \ldots, t_{k-\ell}\}$ that avoids edges with weight zero, and that uses $p_m$, must have exactly one edge with negative weight (there are no edges with negative weight above or to the left of the the $[m,k]$ edge).

Suppose we have a path $p'_m$ from $s_m$ to $t_k$ that avoids edges with weight zero and that makes its first vertical turn before the $[m,k_0]$ edge, i.e,  before column $k_0$, say at column $k_0'<k_0$. Necessarily, $k_0' \leq k_1$. So any path from $s_{m-1}$ to $t_{k-1}$ that avoids edges with weight zero and that is disjoint from $p_m'$ must make its first vertical turn earlier than $p_{m-1}$ does. Let $p_{m-1}'$ be such a path, turning vertical first at column $k_1'<k_1$. We have $k_1' \leq k_2$, so any path from $s_{m-2}$ to $t_{k-2}$ that avoids edges with weight zero and that is disjoint from $p_{m-1}'$ must make its first vertical turn earlier than $p_{m-2}$ does. Iterating, we eventually get that if $p_m', \ldots, p_{m-\ell+1}'$ is any vertex disjoint path system from $\{s_m, \ldots, s_{m-\ell+1}\}$ to $\{t_k, \ldots, t_{k-\ell+1}\}$ that avoids edges with weight zero, then it must be the case that $p_{m-\ell+1}'$ makes its first vertical turn earlier than $p_{m-\ell+1}$ does. So any path from $s_{k-\ell}$ to $t_{k-\ell}$ that is disjoint from $p_{m-\ell+1}'$ and avoids edges with weight zero must make its first vertical turn before it reaches the $[m-\ell, k_\ell]$ edge. But no such vertical turn exists --- if it did, the process of marking edges would not have stopped with the $[m-\ell,k_\ell]$ edge. This contradiction finishes the proof of the assertion that all the vertex disjoint path systems from $\{s_m, \ldots, s_{m-\ell}\}$ to $\{t_k, \ldots, t_{k-\ell}\}$, that avoid edges with weight zero, use the $[m,k_0]$ edge and no other negative edges, and finishes the verification of item 4. 

\subsection{Deriving Theorem \ref{thm-mainGZ} from Claim \ref{claim-alg-works}} \label{subsec-main-proof}

First suppose that ${\bf e}$ is a restricted growth sequence relative to ${\bf a}$. Then it is easy to see that Algorithm \ref{alg-main} terminates with the report that $M_{{\bf e}\rightarrow {\bf a}}$ is TNN, and moreover that we can recover the final weighted planar network $A$ of Algorithm \ref{alg-main}, by the following process: for $i=1, \ldots, n$, in the $i$th iteration of Algorithm \ref{alg-res-growth},
\begin{itemize}
\item if no $e' \in X$ is found with $e' = a_i$ (so all $e' < a_i$), then for each $m=i,\ldots, n$, at the $[m,i]$ edge of the weighted planar network put the weight $a_i-(X)_{m-i+1}$. Here $(X)_j$ indicates the $j$th entry of $X$ (in its $i$th incarnation); and
\item if there is an $e' \in X$ with $e' = a_i$, say $e'=(X)_{m'-i+1}$ for some $m'$, then for each $m=i,\ldots, m'$, at the $[m,i]$ edge of the weighted planar network put the weight $a_i-(X)_{m-i+1}$, and for each $m=m'+1, \ldots, n$, at the $[m,i]$ edge put the weight $0$.  
\end{itemize}
Note moreover that using this process, the weighted planar network that witnesses that $M_{{\bf e}\rightarrow {\bf a}}$ is TNN can be constructed from ${\bf a}$ and ${\bf e}$ in only $O(n^2)$ operations.

On the other hand, suppose that ${\bf e}$ is not a restricted growth sequence relative to ${\bf a}$. If we execute the procedure described above (that associates a weighted planar network with a running of Algorithm \ref{alg-res-growth}), stopping it at the point where an $e' > a_i$ is first found, then the resulting partial weighted planar network agrees with the final network of Algorithm \ref{alg-main} at all the edges that are involved in the collections of vertex disjoint paths that witness a negative minor of $M_{{\bf e}\rightarrow {\bf a}}$, and so the collection of sources and sinks (rows and columns of $M_{{\bf e}\rightarrow {\bf a}}$) that witnesses a negative minor of $M_{{\bf e}\rightarrow {\bf a}}$ can be recovered from the output of Algorithm \ref{alg-res-growth}. 

\section*{Acknowledgement}

We thank an anonymous referee for a very careful reading and many helpful suggestions to improve the clarity of the presentation.

\end{document}